\documentclass[11pt,a4paper]{amsart}
\usepackage[all]{xy}

\newcommand{\ie}{{\em i.e.}\ }
\newcommand{\cf}{{\em cf.}\ }
\newcommand{\eg}{{\em e.g.}\ }
\newcommand{\ko}{\: , \;}
\newcommand{\ul}[1]{\underline{#1}}
\newcommand{\ol}[1]{\overline{#1}}
\renewcommand{\tilde}[1]{\widetilde{#1}}

\newtheorem*{theorem}{Theorem}
\newtheorem*{lemma}{Lemma}
\newtheorem*{proposition}{Proposition}
\newtheorem*{corollary}{Corollary}

\newcommand{\opname}[1]{\operatorname{\mathsf{#1}}}

\renewcommand{\mod}{\opname{mod}\nolimits}
\newcommand{\grmod}{\opname{grmod}\nolimits}
\newcommand{\proj}{\opname{proj}\nolimits}
\newcommand{\Mod}{\opname{Mod}\nolimits}

\newcommand{\per}{\opname{per}\nolimits}

\newcommand{\op}{^{op}}
\newcommand{\der}{\cd}

\newcommand{\dgcat}{\opname{dgcat}\nolimits}
\newcommand{\Ho}{\opname{Ho}\nolimits}
\newcommand{\rep}{\opname{rep}\nolimits}
\newcommand{\enh}{\opname{enh}\nolimits}
\newcommand{\bim}{\opname{rep}\nolimits}
\newcommand{\eff}{\opname{eff}\nolimits}

\newcommand{\pretr}{\opname{pretr}\nolimits}

\newcommand{\ac}{\mathcal{A}c}

\newcommand{\colim}{\opname{colim}}
\newcommand{\mylim}{\opname{lim}}

\newcommand{\Z}{\mathbf{Z}}
\newcommand{\N}{\mathbf{N}}
\newcommand{\ra}{\rightarrow}
\newcommand{\iso}{\stackrel{_\sim}{\rightarrow}}
\newcommand{\id}{\mathbf{1}}

\newcommand{\Hom}{\opname{Hom}}
\newcommand{\RHom}{\opname{RHom}}

\newcommand{\Ext}{\opname{Ext}}

\newcommand{\R}{\mathbf{R}}
\newcommand{\ten}{\otimes}
\newcommand{\lten}{\overset{\boldmath{L}}{\ten}}

\newcommand{\ca}{{\mathcal A}}
\newcommand{\cb}{{\mathcal B}}
\newcommand{\cc}{{\mathcal C}}
\newcommand{\cd}{{\mathcal D}}
\newcommand{\ce}{{\mathcal E}}
\newcommand{\cF}{{\mathcal F}}
\newcommand{\cg}{{\mathcal G}}
\newcommand{\ch}{{\mathcal H}}
\newcommand{\ci}{{\mathcal I}}

\newcommand{\cm}{{\mathcal M}}
\newcommand{\cn}{{\mathcal N}}
\newcommand{\co}{{\mathcal O}}
\newcommand{\cp}{{\mathcal P}}

\newcommand{\ct}{{\mathcal T}}
\newcommand{\cu}{{\mathcal U}}
\newcommand{\cv}{{\mathcal V}}

\newcommand{\eps}{\varepsilon}

\begin{document}

\title[On triangulated orbit categories]{On triangulated orbit categories}

\author{Bernhard Keller}
\address{UFR de Math\'ematiques\\
   UMR 7586 du CNRS \\
   Case 7012\\
   Universit\'e Paris 7\\
   2, place Jussieu\\
   75251 Paris Cedex 05\\
   France }

\email{
\begin{minipage}[t]{5cm}
keller@math.jussieu.fr \\
www.math.jussieu.fr/ $\tilde{ }$ keller
\end{minipage}
}

 \date{June 9, 2003, last
modified on December 5, 2005} 

\dedicatory{Dedicated to Claus Michael Ringel on the occasion of his
sixtieth birthday}

\begin{abstract}
We show that the category of orbits of the bounded derived
category of a hereditary category under a well-behaved
autoequivalence is canonically triangulated. This answers a question
by Aslak Buan, Robert Marsh and Idun Reiten which appeared in
their study \cite{BuanMarshReinekeReitenTodorov04}
with M.~Reineke
and G.~Todorov of the link between tilting theory and cluster
algebras (\cf also \cite{CalderoChapotonSchiffler04})
and a question by Hideto Asashiba about orbit
categories. We observe that the resulting triangulated orbit categories
provide many examples of triangulated categories with the
Calabi-Yau property. These include the category of projective modules
over a preprojective algebra of generalized Dynkin type in the sense
of Happel-Preiser-Ringel \cite{HappelPreiserRingel80},
whose triangulated structure goes back to Auslander-Reiten's
work \cite{AuslanderReiten87}, \cite{Reiten87}, \cite{AuslanderReiten96}.
\end{abstract}

\maketitle

\section{Introduction}
\label{s:introduction} Let $\ct$ be an additive category and
$F:\ct\ra\ct$ an automorphism (a standard construction allows one
to replace a category with autoequivalence by a category with
automorphism). Let $F^\Z$ denote the group of
automorphisms generated by $F$. By definition, the {\em orbit
category} $\ct/F=\ct/F^\Z$ has the same objects as $\ct$ and its
morphisms from $X$ to $Y$ are in bijection with
\[
\bigoplus_{n\in\Z} \Hom_{\ct}(X, F^n Y).
\]
The composition is defined in the natural way
(\cf \cite{CibilsMarcos03}, where this category is called
the skew category). The
canonical projection functor $\pi: \ct \to \ct/F$
is endowed with a natural isomorphism $\pi\circ F \iso \pi$
and 2-universal among such functors.
Clearly $\ct/F$ is still an additive category and the projection is an
additive functor. Now suppose that $\ct$ is a triangulated
category and that $F$ is a triangle functor. Is there a
triangulated structure on the orbit category such that the
projection functor becomes a triangle functor? One can
show that in general, the answer is negative. A closer look at the
situation even gives one the impression that quite strong assumptions are
needed for the answer to be positive. In this article, we give a
sufficient set of conditions. Although they are very strong, they
are satisfied in certain cases of interest. In particular,
one obtains that the cluster categories of
\cite{BuanMarshReinekeReitenTodorov04}, \cite{CalderoChapotonSchiffler04}
are triangulated. One also obtains that the category
of projective modules over the preprojective algebra of a
generalized Dynkin diagram in the sense of
Happel-Preiser-Ringel \cite{HappelPreiserRingel80} is triangulated,
which is also immediate from Auslander-Reiten's
work \cite{AuslanderReiten87}, \cite{Reiten87}, \cite{AuslanderReiten96}.
More generally, our method yields many easily constructed
examples of triangulated categories with the Calabi-Yau property.

Our proof consists in constructing, under quite general
hypotheses, a `triangulated hull' into which the orbit
category $\ct/F$ embeds. Then we show that under
certain restrictive assumptions, its image in the triangulated hull
is stable under extensions and hence equivalent to
the hull.

The contents of the article are as follows: In
section~\ref{s:examples}, we show by examples that triangulated
structures do not descend to orbit categories in general. In
section~\ref{s:main-theorem}, we state the main theorem for
triangulated orbit categories of derived categories of hereditary
algebras. We give a first, abstract, construction of the triangulated
hull of an orbit category in section~\ref{s:triangulated-hull}. This
construction is based on the formalism of dg categories as developped in
\cite{Keller94} \cite{Drinfeld04} \cite{Toen04}. Using the natural
$t$-structure on the derived category of a hereditary category
we prove the main theorem in
section~\ref{s:closure-under-extensions}.

We give a more concrete
construction of the triangulated hull of the orbit category in
section~\ref{s:another-construction}. In some sense, the second
construction is `Koszul-dual' to the first: whereas the first
construction is based on the tensor algebra
\[
T_A(X) = \oplus_{n=0}^\infty X^{\ten_A n}
\]
of a (cofibrant) differential graded bimodule $X$ over a
differential graded algebra $A$, the second one uses the
`exterior algebra'
\[
A \oplus X^\wedge[-1]
\]
on its dual $X^\wedge=\RHom_A(X_A,A)$ shifted by one degree.
In the cases considered by Buan
et al.~\cite{BuanMarshReinekeReitenTodorov04} and
Caldero-Chapoton-Schiffler~\cite{CalderoChapotonSchiffler04}, this
also yields an interesting new description of the orbit category
itself in terms of the stable category \cite{Krause05} of a
differential graded algebra.

In section~\ref{s:On-the-Calabi-Yau-property}, we observe that
triangulated orbit categories provide easily constructed examples of
triangulated categories with the Calabi-Yau property.  Finally, in
section~\ref{s:Universal-properties}, we characterize our
constructions by universal properties in the $2$-category of enhanced
triangulated categories. This also allows us to examine their
functoriality properties and to formulate a more
general version of the main theorem which applies to derived
categories of hereditary categories which are not necessarily module
categories.

\section{Acknowledgments}
The author thanks C.~M.~Ringel for his interest and encouragement
throughout the years.  He is indebted to A.~Buan, R.~Marsh and
I.~Reiten for the question which motivated the main theorem of this
article and to H.~Asashiba for suggesting a first generalization of
its original formulation. He is grateful to H.~Lenzing for pointing
out the need to weaken its hypotheses to include orbit categories of
general hereditary categories. He expresses his thanks to C.~Geiss and
J.-A. de la Pe\~na for the hospitality he enjoyed at the Instituto de
Matem\'aticas, UNAM, in August 2004 while the material of
section~\ref{s:On-the-Calabi-Yau-property} was worked out.  He thanks
R.~Bocklandt, W.~Crawley-Boevey, K.~Erdmann, C.~Geiss, A.~Neeman,
I.~Reiten, A.~Rudakov, J.~Schr\"oer, G.~Tabuada,
M.~Van den Bergh, Feng Ye, Bin Zhu and an anonymous referee
for helpful comments and corrections.

\section{Examples}
\label{s:examples}
Let $\ct$ be a triangulated category, $F:\ct\to\ct$ an autoequivalence
and $\pi: \ct\to\ct/F$ the projection functor. In general, a morphism
\[
u: X \to Y
\]
of $\ct/F$ is given by a morphism
\[
X \to \bigoplus_{i=1}^N F^{n_i} Y
\]
with $N$ non vanishing components $u_1, \ldots, u_N$ in $\ct$. Therefore,
in general, $u$ does not lift to a morphism in $\ct$ and it is not
obvious how to construct a \lq triangle\rq
\[
\xymatrix{X \ar[r]^{u} & Y \ar[r] & Z \ar[r] & SZ}
\]
in $\ct/F$. Thus, the orbit category $\ct/F$ is certainly not
trivially triangulated. Worse, in \lq most\rq\ cases, it is impossible
to endow $\ct/F$ with a triangulated structure such that the
projection functor becomes a triangle functor. Let us consider
three examples where $\ct$ is the bounded derived category
$\der^b(A)=\der^b(\mod A)$ of the category of finitely generated (right) modules
$\mod A$ over an algebra $A$ of finite dimension over a field $k$.
Thus the objects of $\der^b(A)$ are the complexes
\[
M=(\ldots \to M^p \to M^{p+1} \to \ldots)
\]
of finite-dimensional $A$-modules such that $M^p=0$ for all large $|p|$ and morphisms
are obtained from morphisms of complexes by formally inverting all
quasi-isomorphisms. The suspension functor $S$ is defined by $SM=M[1]$, where $M[1]^p=M^{p+1}$
and $d_{M[1]}=-d_M$, and the triangles are constructed from short
exact sequences of complexes.

Suppose that $A$ is hereditary. Then the orbit category $\der^b(A)/S^2$,
first introduced by D.~Happel in \cite{Happel87}, is triangulated.
This result is due to Peng and Xiao \cite{PengXiao97}, who show
that the orbit category is equivalent to the homotopy category
of the category of $2$-periodic complexes of projective $A$-modules.

On the other hand, suppose that $A$ is the algebra of dual
numbers $k[X]/(X^2)$. Then the orbit category $\der^b(A)/S^2$
is not triangulated. This is an observation due to A.~Neeman
(unpublished). Indeed, the
endomorphism ring of the trivial module $k$ in the orbit
category is isomorphic to a polynomial ring $k[u]$.
One checks that the endomorphism $1+u$ is monomorphic.
However, it does not admit a left inverse (or else it
would be invertible in $k[u]$). But in a triangulated
category, each monomorphism admits a left inverse.

One might think that this phenomenon is linked to the
fact that the algebra of dual numbers is of infinite global
dimension. However, it may also occur for algebras of
finite global dimension: Let $A$ be such an algebra.
Then, as shown by D.~Happel in \cite{Happel87}, the
derived category $\der^b(A)$ has Auslander-Reiten triangles.
Thus, it admits an autoequivalence, the
Auslander-Reiten translation $\tau$, defined by
\[
\Hom(?,S\tau M) \iso D\Hom(M,?) \ko
\]
where $D$ denotes the functor $\Hom_k(?,k)$.  Now let $Q$ be
the Kronecker quiver
\[
\xymatrix{1 \ar@/^/[r] \ar@/_/[r] & 2 } .
\]
The path algebra $A=kQ$ is finite-dimensional and hereditary
so Happel's theorem applies. The endomorphism ring of the
image of the free module $A_A$ in the orbit category
$\cd^b(A)/\tau$ is the preprojective algebra $\Lambda(Q)$
(\cf section~\ref{ss:Projectives-over-the-preprojective-algebra}).
Since $Q$ is not a Dynkin quiver, it is infinite-dimensional
and in fact contains a polynomial algebra (generated by
any non zero morphism from the simple projective $P_1$ to 
$\tau^{-1}P_1$). As above, it follows that the
orbit category does not admit a triangulated structure.

\section{The main theorem}
\label{s:main-theorem}

Assume that $k$ is a field, and $\ct$ is the bounded derived
category $\der^b(\mod A)$ of the category of finite-dimensional
(right) modules $\mod A$ over a finite-dimensional $k$-algebra $A$.
Assume that $F:\ct\ra\ct$ is a standard equivalence \cite{Rickard91},
\ie $F$ is isomorphic to the derived tensor product
\[
?\lten_A X : \der^b(\mod A) \to \der^b(\mod A)
\]
for some complex $X$ of $A$-$A$-bimodules.
All autoequivalences with an \lq algebraic construction\rq\
are of this form, \cf section~\ref{s:Universal-properties}.

\begin{theorem} Assume that the following
hypotheses hold:
\begin{itemize}
\item[1)] There is a hereditary abelian $k$-category $\ch$ and a
triangle equivalence
\[
\der^b(\mod A) \iso \der^b(\ch)
\]
In the conditions 2) and 3) below, we identify $\ct$ with $\der^b(\ch)$.
\item[2)] For each indecomposable $U$ of $\ch$, only finitely many
objects $F^i U$, $i\in\Z$, lie in $\ch$.
\item[3)] There is an integer $N\geq 0$ such that the $F$-orbit of
each indecomposable of $\ct$ contains an object $S^n U$,
for some $0\leq n\leq N$ and some indecomposable object $U$ of
$\ch$.
\end{itemize}
Then the orbit category $\ct/F$ admits a natural triangulated
structure such that the projection functor $\ct\to\ct/F$ is
triangulated.
\end{theorem}

The triangulated structure on the orbit category is
(most probably) not unique. However, as we will see in
section~\ref{ss:exact-dg-orbit-categories}, the orbit category is the
associated triangulated category of a dg category, the {\em
exact (or pretriangulated) dg orbit category}, and the exact dg orbit
category is unique and functorial (in the homotopy category of dg categories)
since it is the solution of a universal problem. Thus, although
perhaps not unique, the triangulated structure on the orbit category is
at least canonical, insofar as it comes from a dg structure
which is unique up to quasi-equivalence.

The construction of the triangulated orbit category $\ct/F$ via the exact
dg orbit category also shows that there is a triangle equivalence
between $\ct/F$ and the stable category $\ul{\ce}$ of some
Frobenius category $\ce$.

In sections~\ref{ss:The-motivating-example},
\ref{ss:Projectives-over-the-preprojective-algebra} and
\ref{ss:Projectives-over-Ln}, we will illustrate the
theorem by examples. In sections~\ref{s:triangulated-hull} and
\ref{s:closure-under-extensions} below, we prove the theorem. The
strategy is as follows: First, under very weak assumptions,
we embed $\ct/F$ in a naturally
triangulated ambient category $\cm$ (whose intrinsic
interpretation will be given in section~\ref{ss:exact-dg-orbit-categories}).
Then we show that $\ct/F$ is closed under extensions in 
the ambient category $\cm$.
Here we will need the full strength of the assumptions 1), 2) and 3).

If $\ct$ is the derived category of an abelian category which
is not necessarily a module category, one can still
define a suitable notion of a standard equivalence $\ct\to\ct$,
\cf section~\ref{s:Universal-properties}. Then the
analogue of the above theorem is true, \cf
section~\ref{ss:hereditary-categories}.

\section{Construction of the triangulated hull $\cm$}
\label{s:triangulated-hull}

The construction is based on the formalism of dg categories,
which is briefly recalled in section~\ref{ss:htpy-cat-dgcat}. We refer to
\cite{Keller94}, \cite{Drinfeld04} and \cite{Toen04} for more
background.

\subsection{The dg orbit category} \label{ss:dgorbitcategory}
Let $\ca$ be a dg category and $F: \ca \ra \ca$ a
dg functor inducing an equivalence in $H^0 \ca$. 
We define the {\em dg orbit category $\cb$} to be the dg
category with the same objects as $\ca$ and such that for $X,Y\in
\cb$, we have
\[
\cb\,(X,Y) \iso \colim_p \bigoplus_{n\geq 0} \ca\,(F^n X,
F^{p} Y) \ko
\]
where the transition maps of the $\colim$ are given by $F$.
This definition ensures that $H^0\cb$ is isomorphic to the orbit
category $(H^0\ca)/F$.

\subsection{The projection functor and its right adjoint}
\label{ss:projection-functor}
From now on, we assume that
for all objects $X,Y$ of $\ca$, the group
\[
(H^0\ca)\,(X,F^n Y)
\]
vanishes except for finitely many $n\in\Z$.
We have a canonical dg functor $\pi : \ca \ra \cb$. It yields an
$\ca$-$\cb$-bimodule
\[
(X,Y) \mapsto \cb\,(\pi X, Y).
\]
The standard functors associated with this bimodule are
\begin{itemize}
\item[-] the derived tensor functor (=induction functor)
\[
\pi_* : \der\ca \ra \der\cb
\]
\item[-] the derived Hom-functor (right adjoint to $\pi_*$), which
equals the restriction along~$\pi$:
\[
\pi_\rho : \der\cb \ra \der\ca
\]
\end{itemize}
For $X\in \ca$, we have
\[
\pi_* (X^\wedge) = (\pi X)^\wedge \ko
\]
where $X^\wedge$ is the functor represented by $X$.
Moreover, we have an isomorphism in $\cd\ca$
\[
\pi_\rho \pi_* (X^\wedge) = \bigoplus_{n\in\Z} F^n (X^\wedge),
\]
by the definition of the morphisms of $\cb$ and the vanishing
assumption made above.

\subsection{Identifying objects of the orbit category} \label{ss:identifying}
The functor $\pi_\rho : \der\cb \ra\der\ca$ is the restriction along
a morphism of dg categories. Therefore, it detects isomorphisms.
In particular, we obtain the following: Let $E\in\der\cb$,
$Z\in\der\ca$ and let $f : Z \ra \pi_\rho E$ be a morphism. Let $g:
\pi_* Z \ra E$ be the morphism corresponding to $f$ by the
adjunction.  In order to show that $g$ is an isomorphism, it is
enough to show that $\pi_\rho g: \pi_\rho \pi_* Z \ra \pi_\rho E$
is an isomorphism.

\subsection{The ambient triangulated category}
\label{ss:ambient}
We use the notations of the main theorem. Let $X$ be a complex of
$A$-$A$-bimodules such that $F$ is isomorphic to the total derived
tensor product by $X$. We may assume that $X$ is bounded and that its
components are projective on both sides. Let $\ca$ be the dg category
of bounded complexes of finitely generated projective $A$-modules. The
tensor product by $X$ defines a dg functor from $\ca$ to $\ca$. By
abuse of notation, we denote this dg functor by $F$ as well. The
assumption 2) implies that the vanishing assumption of
subsection~\ref{ss:projection-functor} is satisfied. Thus we obtain a dg
category $\cb$ and an equivalence of categories
\[
\der^b(\mod A)/ F \iso H^0 \cb.
\]
We let the ambient triangulated category
$\cm$ be the triangulated subcategory of $\der\cb$ generated
by the representable functors. The Yoneda embedding $H^0 \cb \ra \der\cb$
yields the canonical embedding $\der^b(\mod A)/F \ra \cm$. We have
a canonical equivalence $\der(\Mod A) \iso \der\ca$ and therefore we
obtain a pair of adjoint functors $(\pi_*, \pi_\rho)$ between
$\der(\Mod A)$ and $\cd\cb$. The functor $\pi_*$ restricts to the
canonical projection $\der^b(\mod A) \ra \der^b(\mod A)/F$.

\section{The orbit category is closed under extensions}
\label{s:closure-under-extensions}

Consider the right adjoint $\pi_\rho$ of $\pi_*$. It is defined on
$\der\cb$ and takes values in the unbounded derived category of
all $A$-modules $\der(\Mod A)$. For $X \in \der^b(\mod A)$, the object
$\pi_* \pi_\rho X$ is isomorphic to the sum of the translates $F^i
X$, $i\in\Z$, of $X$. It follows from assumption 2) that for each
fixed $n\in\Z$, the module $H^n_{\mod A} F^i X$ vanishes for
almost all $i\in\Z$. Therefore the sum of the $F^i X$ lies in
$\der(\mod A)$.

Consider a morphism $f: \pi_* X \ra \pi_* Y$ of the orbit category
$\ct/F=\der^b(\mod A)/F$. We form a triangle
\[
\pi_* X \ra \pi_* Y \ra E \ra S\pi_* X
\]
in $\cm$. We apply the right adjoint $\pi_\rho$ of $\pi_*$ to this
triangle. We get a triangle
\[
\pi_\rho \pi_* X \ra \pi_\rho \pi_* Y \ra \pi_\rho E \ra S
\pi_\rho \pi_* X
\]
in $\der(\Mod A)$. As we have just seen, the terms $\pi_\rho \pi_* X$ and
$\pi_\rho \pi_* Y$ of the triangle belong to $\der(\mod A)$.  Hence so
does $\pi_\rho E$. We will construct an object $Z\in \der^b(\mod A)$ and
an isomorphism $g: \pi_* Z \ra E$ in the orbit category. Using
proposition~\ref{ss:Extension} below, we extend the canonical
$t$-structure on $\der^b(\ch) \iso \der^b(\mod A)$ to a $t$-structure on
all of $\der(\mod A)$. Since $\ch$ is hereditary, part b) of the
proposition shows that each object of $\der(\mod A)$ is the sum of its
$\ch$-homology objects placed in their respective degrees. In
particular, this holds for $\pi_\rho E$. Each of the homology objects
is a finite sum of indecomposables.  Thus $\pi_\rho E$ is a sum of
shifted copies of indecomposable objects of $\ch$. Moreover, $F
\pi_\rho E$ is isomorphic to $\pi_\rho E$, so that the sum is stable
under $F$. Hence it is a sum of $F$-orbits of shifted indecomposable
objects. By assumption 3), each of these orbits contains an
indecomposable direct factor of
\[
\bigoplus_{0\leq n \leq N} (H^n \pi_\rho E)[-n].
\]
Thus there are only finitely many orbits involved. Let
$Z_1,\ldots,Z_M$ be shifted indecomposables of $\ch$ such that
$\pi_\rho E$ is the sum of the $F$-orbits of the $Z_i$. Let $f$ be
the inclusion morphism
\[
Z=\bigoplus_{i=1}^M Z_i \ra \pi_\rho E
\]
and $g : \pi_* Z \ra E$ the morphism corresponding to $f$ under
the adjunction. Clearly $\pi_\rho g$ is an isomorphism. By
subsection~\ref{ss:identifying}, the morphism $g$ is invertible
and we are done.

\subsection{Extension of $t$-structures to unbounded categories}
\label{ss:Extension}
Let $\ct$ be a triangulated category and $\cu$ an aisle \cite{KellerVossieck88}
in $\ct$. Denote the associated $t$-structure \cite{BeilinsonBernsteinDeligne82}
by $(\cu_{\leq 0}, \cu_{\geq 0})$,
its heart by $\cu_0$, its homology functors by $H^n_\cu : \ct \to
\cu_0$ and its truncation functors by $\tau_{\leq n}$ and $\tau_{>n}$.
Suppose that $\cu$ is {\em dominant}, \ie the following two
conditions hold:
\begin{itemize}
\item[1)] a morphism $s$ of $\ct$ is
an isomorphism iff $H^n_\cu(s)$ is an isomorphism for all $n\in\Z$ and
\item[2)] for each object $X\in\ct$, the canonical morphisms
\[
\Hom(?,X) \to \mylim\Hom(?, \tau_{\leq n} X)
\mbox{ and }
\Hom(X,?) \to \mylim\Hom(\tau_{>n} X, ?)
\]
are surjective.
\end{itemize}
Let $\ct^b$ be the full triangulated subcategory of
$\ct$ whose objects are the $X\in\ct$ such that $H^n_\cu(X)$ vanishes
for all $|n|\gg 0$. Let $\cv^b$ be an aisle on $\ct^b$. Denote
the associated $t$-structure on $\ct^b$ by $(\cv_{\leq n}, \cv_{>n})$,
its heart by $\cv_0$, the homology functor by $H^n_{\cv^b} : \cb \to \cv_0$
and its truncation functors by $(\sigma_{\leq 0}, \sigma_{>0})$.

Assume that there is an $N\gg 0$ such that we have
\[
H^0_{\cv^b} \iso H^0_{\cv^b} \tau_{>-n} \mbox{ and }
H^0_{\cv^b} \tau_{\leq n} \iso H^0_{\cv^b}
\]
for all $n\geq N$. We define $H^0_\cv: \ct \to \cv_0$ by
\[
H^0_\cv(X) = \colim H^0_{\cv^b} \tau_{>-n} \tau_{\leq m} X
\]
and $H^n_\cv(X)= H^0_\cv S^n X$, $n\in\Z$.
We define $\cv\subset \ct$ to be the full subcategory of $\ct$
whose objects are the $X\in\ct$ such that $H^n_\cv(X)=0$ for all $n>0$.

\begin{proposition}
\begin{itemize} \item[a)] $\cv$ is an aisle in $\ct$ and the associated
$t$-structure is dominant.
\item[b)] If $\cv^b$ is {\em hereditary}, \ie each triangle
\[
\sigma_{\leq 0} X \to X \to \sigma_{>0} X \to S \sigma_{\leq 0} X \ko X\in\ct^b
\]
splits, then $\cv$ is hereditary and each object $X\in\ct$ is (non canonically)
isomorphic to the sum of the $S^{-n}H^n_\cv(X)$, $n\in\Z$.
\end{itemize}
\end{proposition}
The proof is an exercise on $t$-structures which we leave to the reader.

\section{Another construction of the triangulated hull of the orbit category}
\label{s:another-construction}

\subsection{The construction} Let $A$ be a finite-dimensional algebra of finite global
dimension over a field $k$. Let $X$ be an $A$-$A$-bimodule complex
whose homology has finite total dimension.
Let $F$ be the functor
\[
? \lten_A X : \der^b(\mod A) \to \der^b(\mod A).
\]
We suppose that $F$ is an equivalence and that
for all $L,M$ in $\der^b(\mod A)$, the group
\[
\Hom(L,F^n M)
\]
vanishes for all but finitely many $n\in\Z$.
We will construct a triangulated category equivalent
to the triangulated hull of section~\ref{s:triangulated-hull}.

Consider $A$ as a dg algebra concentrated in degree $0$.
Let $B$ be the dg algebra with underlying complex
$A \oplus X[-1]$, where the multiplication is that of the
trivial extension:
\[
(a,x)(a',x') = (aa', ax'+xa').
\]
Let $\der B$ be the derived category of $B$ and $\der^b(B)$ the {\em
bounded derived category}, \ie the full subcategory of $\der B$ formed
by the dg modules whose homology has finite total dimension over
$k$. Let $\per(B)$ be the {\em perfect} derived category of $B$, \ie
the smallest subcategory of $\der B$ containing $B$ and stable under
shifts, extensions and passage to direct factors.  By our assumption
on $A$ and $X$, the perfect derived category is contained in $\der^b(B)$.
The obvious morphism $B \to A$ induces a restriction functor
$\der^b A \to \der^b B$ and by composition, we obtain a functor
\[
\der^b A \to \der^b B \to \der^b(B)/ \per(B)
\]

\begin{theorem}
The category $\der^b(B)/\per(B)$ is equivalent to the
triangulated hull (\cf section~\ref{s:triangulated-hull})
of the orbit category of $\der^b(A)$ under $F$ and
the above functor identifies with the projection functor.
\end{theorem}

\begin{proof} If we replace $X$ by a quasi-isomorphic bimodule,
the algebra $B$ is replaced by a quasi-isomorphic dg algebra
and its derived category by an equivalent one. Therefore,
it is no restriction of generality to assume, as we will
do, that $X$ is cofibrant as a dg $A$-$A$-bimodule. We will
first compute morphisms in $\der B$ between dg $B$-modules
whose restrictions to $A$ are cofibrant. For this, let
$C$ be the dg submodule of the bar resolution of $B$ as a bimodule
over itself whose underlying graded module is
\[
C= \coprod_{n\geq 0} B\ten_A X^{\ten_A n}\ten_A B.
\]
The bar resolution of $B$ is a coalgebra in the category
of dg $B$-$B$-bimodules (\cf \eg \cite{Keller94}) and $C$ becomes
a dg subcoalgebra. Its counit is
\[
\eps : C \to B\ten_A B \to B
\]
and its comultiplication is given by
\[
\Delta(b_0, x_1, \ldots, x_n, b_{n+1}) = \sum_{i=0}^n
(b_0, x_1, \ldots, x_i) \ten 1\ten 1 \ten (x_{i+1}, \ldots, b_{n+1}).
\]
It is not hard to see that the inclusion of $C$ in the bar
resolution is an homotopy equivalence of left (and of right)
dg $B$-modules. Therefore, the same holds for the counit $\eps : C \to B$.
For an arbitrary right dg $B$-module $L$, the counit $\eps$ thus
induces a quasi-isomorphism $L\ten_B C \to L$. Now suppose that
the restriction of $L$ to $A$ is cofibrant. Then $L\ten_A B \ten_A B$
is cofibrant over $B$ and thus $L\ten_B C \to L$ is a cofibrant
resolution of $L$. Let $\cc_1$ be the dg category whose objects are
the dg $B$-modules whose restriction to $A$ is cofibrant and
whose morphism spaces are the
\[
\Hom_B(L\ten_B C, M\ten_B C).
\]
Let $\cc_2$ be the dg category with the same objects as $\cc_1$ and
whose morphism spaces are
\[
\Hom_B(L\ten_B C, M).
\]
By definition, the composition of two morphisms $f$ and $g$ of $\cc_2$ is
given by
\[
f\circ (g\ten \id_C) \circ (\id_L \ten \Delta).
\]
We have a dg functor $\Phi: \cc_2 \to \cc_1$ which is the identity on objects
and sends $g: L \to M$ to
\[
(g\ten \id_C) \circ (\id_L \ten \Delta) : L\ten_B C \to M\ten_B C.
\]
The morphism
\[
\Hom_B(L\ten_B C, M) \to \Hom_B(L\ten_B C, M\ten_B C)
\]
given by $\Phi$ is left inverse to the quasi-isomorphism induced
by $M\ten_B C \to M$. Therefore, the dg functor $\Phi$ yields a
quasi-isomorphism between $\cc_2$ and $\cc_1$ so that we can
compute morphisms and compositions in $\der B$ using $\cc_2$. Now
suppose that $L$ and $M$ are cofibrant dg $A$-modules. Consider
them as dg $B$-modules via restriction along the projection $B\to A$.
Then we have natural isomorphisms of complexes
\[
\Hom_{\cc_2}(L,M)=\Hom_B(L\ten_B C, M) \iso \Hom_A(\coprod_{n\geq 0}
L\ten_A X^{\ten_A n}, M).
\]
Moreover, the composition of morphisms in $\cc_2$ translates
into the natural composition law for the right hand side.
Now we will compute morphisms in the quotient category
$\der B/ \per(B)$. Let $M$ be as above. For $p\geq 0$, let
$C_{\leq p}$ be the dg subbimodule with underlying graded
module
\[
\coprod_{n=0}^p B \ten_A X^{\ten_A n} \ten_A B.
\]
Then each morphism
\[
P \to M\ten_B C
\]
of $\der B$ from a perfect object $P$ factors through
\[
M\ten_B C_{\leq p} \to M\ten_B C
\]
for some $p\geq 0$. Therefore, the following complex
computes morphisms in $\der B/\per(B)$:
\[
\colim_{p\geq 1} \Hom_B(L\ten_B C, M\ten_B C/C_{\leq p-1}).
\]
Now it is not hard to check that the inclusion
\[
M\ten X^{\ten_A p} \ten_A A \to M \ten_B (C/C_{\leq p-1})
\]
is a quasi-isomorphism of dg $B$-modules. Thus we obtain quasi-isomorphisms
\[
\Hom_B(L\ten_B C, M\ten X^{\ten_A p} \ten_A A) \to \Hom_B(L\ten_B C,
M\ten_B C/C_{\leq p-1})
\]
and
\[
\prod_{n\geq 0} \Hom_A(L\ten_A X^{\ten_A n}, M\ten_A X^{\ten_A p})
\to\Hom_B(L\ten_B C, M\ten_B C/C_{\leq p-1}).
\]
Moreover, it is not hard to check that if we define transition
maps
\[
\prod_{n\geq 0} \Hom_A(L\ten_A X^{\ten_A n}, M\ten_A X^{\ten_A p})
\to
\prod_{n\geq 0} \Hom_A(L\ten_A X^{\ten_A n}, M\ten_A X^{\ten_A (p+1)})
\]
by sending $f$ to $f\ten_A \id_X$, then we obtain a quasi-isomorphism
of direct systems of complexes. Therefore, the following complex
computes morphisms in $\der B/\per(B)$:
\[
\colim_{p\geq 1} \prod_{n\geq 0} \Hom_A(L\ten_A X^{\ten_A n}, M\ten_A X^{\ten_A p}).
\]
Let $\cc_3$ be the dg category whose objects are the cofibrant dg $A$-modules
and whose morphisms are given by the above complexes.
If $L$ and $M$ are cofibrant dg $A$-modules and belong to $\der^b(\mod A)$, then,
by our assumptions on $F$, this complex is quasi-isomorphic to its subcomplex
\[
\colim_{p\geq 1} \coprod_{n\geq 0} \Hom_A(L\ten_A X^{\ten_A n}, M\ten_A X^{\ten_A p}).
\]
Thus we obtain a dg functor
\[
\cb \to \cc_3
\]
(where $\cb$ is the dg category defined in \ref{ss:dgorbitcategory})
which induces a fully faithful functor $H^0(\cb) \to \der B/\per(B)$ and
thus a fully faithful functor $\cm \to \der^b(B)/\per(B)$. This functor
is also essentially surjective. Indeed, every object in $\der^b(B)$ is
an extension of two objects which lie in the image of $\der^b(\mod A)$.
The assertion about the projection functor is clear from the
above proof.
\end{proof}

\subsection{The motivating example}
\label{ss:The-motivating-example}
Let us suppose that the functor $F$ is
given by
\[
M \mapsto \tau S^{-1}M \ko
\]
where $\tau$ is the Auslander-Reiten translation of $\der^b(A)$ and
$S$ the shift functor. This is the case considered in
\cite{BuanMarshReinekeReitenTodorov04} for the construction
of the cluster category. The functor $F^{-1}$ is isomorphic to
\[
M \mapsto S^{-2} \nu M
\]
where $\nu$ is the Nakayama functor
\[
\nu = ?\lten_A DA \ko DA=\Hom_k(A,k).
\]
Thus $F^{-1}$ is given by the bimodule $X=(DA)[-2]$ and $B=A \oplus
(DA)[-3]$ is the trivial extension of $A$ with a non standard
grading: $A$ is in degree $0$ and $DA$
in degree $3$. For example, if $A$ is the quiver algebra
of an alternating quiver whose underlying graph is $A_n$, then
the underlying ungraded algebra of $B$ is the quadratic dual of the
preprojective algebra associated with $A_n$, \cf \cite{BrennerButlerKing02}.
The algebra $B$ viewed as a differential graded algebra
was investigated by Khovanov-Seidel in \cite{KhovanovSeidel02}.
Here the authors show that $\der^b(B)$ admits a canonical
action by the braid group on $n+1$ strings, a result which
was obtained independently in a similar context by Zimmermann-Rouquier
\cite{RouquierZimmermann03}. The canonical generators of
the braid group act by triangle functors $T_i$ endowed with morphisms
$\phi_i : T_i \to \id$. The cone on each $\phi_i$ belongs
to $\per(B)$ and $\per(B)$ is in fact equal to its smallest
triangulated subcategory stable under direct factors and
containing these cones. Thus, the action becomes trivial in
$\der^b(B)/\per(B)$ and in a certain sense, this is the largest quotient where
the $\phi_i$ become invertible.

\subsection{Projectives over the preprojective algebra}
\label{ss:Projectives-over-the-preprojective-algebra}
Let $A$ be the path algebra of a Dynkin quiver, \ie a
quiver whose underlying graph is a Dynkin diagram of type
$A$, $D$ or $E$. Let $C$ be
the associated preprojective algebra \cite{GelfandPonomarev79},
\cite{DlabRingel80}, \cite{Ringel98}.
In Proposition~3.3 of \cite{AuslanderReiten96}, Auslander-Reiten
show that the category of projective modules over $C$ is
equivalent to the stable category of maximal Cohen-Macaulay
modules over a representation-finite isolated hypersurface
singularity. In particular, it is triangulated. This can
also be deduced from our main theorem: Indeed, it follows from
D.~Happel's description~\cite{Happel87} of the derived
category $\der^b(A)$ that the category $\proj C$ of
finite dimensional projective $C$-modules is equivalent
to the orbit category $\der^b(A)/\tau$, \cf also \cite{Geiss04}.
Moreover, by the theorem of the previous section,
we have an equivalence
\[
\proj C \iso \der^b(B)/\per(B)
\]
where $B=A \oplus (DA)[-2]$. This equivalence
yields in fact more than just a triangulated
structure: it shows that $\proj C$ is endowed with
a canonical Hochschild $3$-cocycle $m_3$,
\cf for example \cite{BensonKrauseSchwede04}.
It would be interesting to identify this cocycle in
the description given in \cite{ErdmannSnashall98}.

\subsection{Projectives over $\Lambda(L_n)$}
\label{ss:Projectives-over-Ln}
The category of projective modules over the algebra
$k[\eps]/(\eps^2)$ of dual numbers is triangulated.
Indeed, it is equivalent to the orbit category of the
derived category of the path algebra of a quiver of
type $A_2$ under the Nakayama autoequivalence~$\nu$.
Thus, we obtain examples of triangulated
categories whose Auslander-Reiten quiver contains a loop.
It has been known since Riedtmann's work \cite{Riedtmann80}
that this cannot occur in the stable category
(\cf below) of a selfinjective finite-dimensional algebra.
It may  therefore seem surprising, \cf \cite{XiaoZhu02}, that loops
do occur in this more general context.
However, loops already do occur in stable categories of
finitely generated reflexive
modules over certain non commutative generalizations
of local rings of rational double points, as shown
by Auslander-Reiten in \cite{AuslanderReiten87}.
These were completely classified by Reiten-Van den Bergh
in \cite{ReitenVandenBergh89}.
In particular, the example of the dual numbers and
its generalization below are among the cases covered by
\cite{ReitenVandenBergh89}.

The example of the dual numbers generalizes as follows:
Let $n\geq 1$ be an integer.  Following \cite{HappelPreiserRingel80a},
the {\em generalized Dynkin graph $L_n$} is defined as the graph
\[
\xymatrix{1 \ar@{-}@(ul,dl)[] \ar@{-}[r] & 2 \ar@{-}[r] & &
\cdots & \ar@{-}[r] & n-1 \ar@{-}[r] & n}
\]
Its edges are in natural bijection with the orbits of the
involution which exchanges each arrow $\alpha$ with $\overline{\alpha}$
in the following quiver:
\[
\xymatrix{1 \ar@(ul,dl)[]_{\eps=\ol{\eps}} \ar@<2pt>[r]^{a_1} &
2 \ar@<2pt>[l]^{\ol{a_1}} \ar@<2pt>[r]^{a_2} & \ar@<2pt>[l]^{\ol{a_2}} &
\cdots & \ar@<2pt>[r]^-{a_{n-2}} &
n-1  \ar@<2pt>[l]^-{\ol{a_{n-2}}}
\ar@<2pt>[r]^-{a_{n-1}} & n \ar@<2pt>[l]^-{\ol{a_{n-1}}} } .
\]
The associated preprojective algebra $\Lambda(L_n)$ of generalized
Dynkin type $L_n$
is defined as the quotient of the path algebra of this quiver
by the ideal generated by the relators
\[
r_v=\sum \alpha\,\overline{\alpha} \ko
\]
where, for each $1\leq v\leq n$, the sum ranges over the arrows
$\alpha$ with starting point $v$. Let $A$ be the path
algebra of a Dynkin quiver with underlying Dynkin graph $A_{2n}$.
Using D.~Happel's description~\cite{Happel87} of the
derived category of a Dynkin quiver, we see that
the orbit category $\der^b(A)/(\tau^n S)$ is equivalent to the
category of finitely generated projective modules
over the algebra $\Lambda(L_n)$.
By the main theorem, this category is thus triangulated.
Its Auslander-Reiten quiver is given by the ordinary quiver
of $\Lambda(L_n)$, \cf above, endowed with $\tau=\id$:
Indeed, in $\der^b(A)$, we have $S^2=\tau^{-(2n+1)}$ so that
in the orbit category, we obtain
\[
\id = (\tau^n S)^2 = \tau^{2n} S^2 = \tau^{-1}.
\]

\section{On the Calabi-Yau property}
\label{s:On-the-Calabi-Yau-property}

\subsection{Serre functors and localizations}
Let $k$ be a field and $\ct$ a $k$-linear triangulated
category with finite-dimensional $\Hom$-spaces. We denote
the suspension functor of $\ct$ by $S$.
Recall from \cite{ReitenVandenBergh02}
that a {\em right Serre functor} for $\ct$ is the datum of
a triangle functor $\nu:\ct\to\ct$ together with
bifunctor isomorphisms
\[
D\Hom_\ct(X,?) \iso \Hom_\ct(?,\nu X)\ko X \in \ct\ko
\]
where $D=\Hom_k(?,k)$. If $\nu$ exists, it is unique up to
isomorphism of triangle functors.
Dually, a {\em left Serre functor} is the
datum of a triangle functor $\nu':\ct\to\ct$ and isomorphisms
\[
D\Hom_\ct(?,X) \iso \Hom_\ct(\nu', ?) \ko X \in\ct.
\]
The category $\ct$ {\em has Serre duality} if it has both a left
and a right Serre functor, or equivalently, if it has a onesided
Serre functor which is an equivalence, \cf \cite{ReitenVandenBergh02}
\cite{BondalVandenBergh03}.
The following lemma is used in \cite{BuanMarshReiten04}.

\begin{lemma} Suppose that $\ct$ has a left Serre functor $\nu'$.
Let $\cu\subset\ct$ be a thick triangulated subcategory and
$L : \ct\to\ct/\cu$ the localization functor.
\begin{itemize}
\item[a)] If $L$ admits a right adjoint $R$, then $L\nu' R$ is a
left Serre functor for $\ct/\cu$.
\item[b)] More generally,
if the functor $\nu': \ct\to\ct$ admits
a total right derived functor $\R\nu':\ct/\cu \to \ct/\cu$
in the sense of Deligne \cite{Deligne73} with respect to the
localization $\ct\to\ct/\cu$, then $\R \nu'$ is a left Serre
functor for $\ct/\cu$.
\end{itemize}
\end{lemma}

\begin{proof} a) For $X$, $Y$ in $\ct$, we have
\[
\Hom_{\ct/\cu}(L\nu' R X, Y) =
\Hom_\ct(\nu' R X, RY) =
D\Hom_\ct(RY, RX) =
D\Hom_{\ct/\cu}(Y,X).
\]
Here, for the last isomorphism, we have used that $R$
is fully faithful ($L$ is a localization functor).

b) We assume that $\ct$ is small. Let $\Mod\ct$ denote the (large) category
of functors from $\ct\op$ to the category of abelian groups and
let $h: \ct\to\Mod\ct$ denote the Yoneda embedding. Let $L^*$ be
the unique right exact functor $\Mod \ct\to \Mod (\ct/\cu)$ which
sends $hX$ to $hLX$, $X\in\ct$. By the calculus of (right) fractions,
$L^*$ has a right adjoint $R$ which takes an object $Y$ to
\[
\colim_{\Sigma_Y} hY' \ko
\]
where the $\colim$ ranges over the category $\Sigma_Y$ of
morphisms $s: Y \to Y'$ which
become invertible in $\ct/\cu$. Clearly $L^* R$ is isomorphic to the
identity so that $R$ is fully faithful. By definition of the total
right derived functor, for each object $X\in\ct/\cu$, the functor
\[
\colim_{\Sigma_X} h(L\nu' X') = L^* \nu'^* R h(X)
\]
is represented by $\R \nu'(X)$. Therefore, we have
\[
\Hom_{\ct/\cu}(\R \nu'(X), Y) =
\Hom_{\Mod \ct/\cu} (L^* \nu'^* R h(X), h(Y)) =
\Hom_{\Mod\ct}(\nu'^* R h(X), R h(Y)).
\]
Now by definition, the last term is isomorphic to
\[
\Hom_{\Mod\ct}(\colim_{\Sigma_X} h(\nu' X'), \colim_{\Sigma_Y} hY') =
\mylim_{\Sigma_X}
\colim_{\Sigma_Y}
\Hom_{\ct}(\nu' X', Y')
\]
and this identifies with
\[
\mylim_{\Sigma_X}
\colim_{\Sigma_Y}
D\Hom_\ct(Y',X') =
D(
\colim_{\Sigma_X}
\mylim_{\Sigma_Y}
\Hom_\ct(Y',X')) =
D\Hom_{\ct/\cu}(Y,X).
\]

\end{proof}

\subsection{Definition of the Calabi-Yau property}
Keep the hypotheses of the preceding section.
By definition \cite{Kontsevich98}, the triangulated category
$\ct$ is {\em Calabi-Yau of CY-dimension $d$}
if it has Serre duality and there is an isomorphism
of triangle functors
\[
\nu \iso S^d.
\]
By extension, if we have $\nu^e \iso S^d$ for some integer $e>0$,
one sometimes says that $\ct$ is Calabi-Yau of fractional dimension $d/e$.
Note that $d\in\Z$ is only determined up to a multiple
of the order of $S$. It would be interesting to link the
CY-dimension to Rouquier's \cite{Rouquier03} notion of dimension of a
triangulated category.

The terminology has its origin in the following example:
Let $X$ be a smooth projective variety of dimension $d$ and
let $\omega_X = \Lambda^d T^*_X$ be the canonical bundle.
Let $\ct$ be the bounded derived category of coherent sheaves
on $X$. Then the classical Serre duality
\[
D\Ext^i_X(\cF,\cg) \iso \Ext^{d-i}(\cg, \cF\ten \omega_X) \ko
\]
where $\cF$, $\cg$ are coherent sheaves, lifts to the isomorphism
\[
D\Hom_\ct(\cF, \cg) \iso \Hom_\ct(\cg, \cF\ten\omega_X[d])\ko
\]
where $\cF$, $\cg$ are bounded complexes of coherent sheaves. Thus
$\ct$ has Serre duality and $\nu = ?\ten \omega_X[d]$. So
the category $\ct$ is Calabi-Yau of CY-dimension $d$
iff $\omega_X$ is isomorphic to $\co_X$, which means
precisely that the variety $X$ is Calabi-Yau of dimension $d$.

If $\ct$ is a Calabi-Yau triangulated category \lq of algebraic
origin\rq\ (for example, the derived category of a category of modules
or sheaves), then it often comes from a Calabi-Yau
$A_\infty$-category. These are of considerable interest in
mathematical physics, since, as Kontsevich shows \cite{Kontsevich92},
\cite{Kontsevich2004}, \cf also \cite{Costello04}, a topological
quantum field theory is associated with each Calabi-Yau
$A_\infty$-category satisfying some additional
assumptions\footnote{Namely, the associated triangulated category
should admit a generator whose endomorphism $A_\infty$-algebra $B$ is
compact (\ie finite-dimensional), smooth (\ie $B$ is perfect as a
bimodule over itself), and whose associated Hodge-de Rham spectral
sequence collapses (this property is conjectured to hold for all
smooth compact $A_\infty$-algebras over a field of characteristic
$0$).}.

\subsection{Examples} \label{ss:CY-examples} {\bf (1)} If $A$ is a finite-dimensional
$k$-algebra, then the homotopy category $\ct$ of bounded complexes
of finitely generated projective $A$-modules has a Nakayama
functor iff $DA$ is of finite projective dimension. In
this case, the category $\ct$ has Serre duality iff
moreover $A_A$ is of finite injective dimension,
\ie iff $A$ is Gorenstein, \cf \cite{Happel91}.
Then the category $\ct$ is Calabi-Yau (necessarily
of CY-dimension $0$) iff $A$ is symmetric.

{\bf (2)}
If $\Delta$ is a
Dynkin graph of type $A_n$, $D_n$, $E_6$, $E_7$ or $E_8$
and $h$ is its Coxeter number
(\ie $n+1$, $2(n-1)$, $12$, $18$ or $30$, respectively), then
for the bounded derived category of finitely generated
modules over a quiver with underlying graph $\Delta$,
we have isomorphisms
\[
\nu^h = (S\tau)^h = S^h \tau^h = S^{(h-2)}.
\]
Hence this category is Calabi-Yau of fractional dimension
$(h-2)/h$.

{\bf (3)}
Suppose that $A$ is a finite-dimensional algebra which is
selfinjective (\ie $A_A$ is also an injective $A$-module).  Then the
category $\mod A$ of finite-dimensional $A$-modules is Frobenius, \ie
it is an abelian (or, more generally, an exact) category with enough
projectives, enough injectives and where an object is projective iff
it is injective. The stable category $\ul{\mod}\, A$ obtained by
quotienting $\mod A$ by the ideal of morphisms factoring through
injectives is triangulated, \cf \cite{Happel87}. The inverse of its
suspension functor sends a module $M$ to the kernel $\Omega M$ of an
epimorphism $P \to M$ with projective $P$.  Let $\cn M = M \ten_A
DA$. Then $\ul{\mod}\,A$ has Serre duality with Nakayama functor $\nu
= \Omega \circ \cn$. Thus, the stable category is Calabi-Yau
of CY-dimension $d$ iff we have an isomorphism of triangle functors
\[
\Omega^{(d+1)}\circ \cn = \id.
\]
For this, it is clearly sufficient that we have an isomorphism
\[
\Omega^{d+1}_{A^e}(A)\ten_A DA \iso A \ko
\]
in the stable category of $A$-$A$-bimodules, \ie
modules over the selfinjective algebra $A\ten A\op$.
For example, we deduce that if $A$ is the path
algebra of a cyclic quiver with $n$ vertices divided by the ideal
generated by all paths of length $n-1$, then $\ul{\mod}\, A$ is
Calabi-Yau of CY-dimension $3$.

{\bf (4)} Let $A$ be a dg algebra. Let $\per(A)\subset\der(A)$
be the subcategory of perfect dg $A$-modules, \ie the
smallest full triangulated subcategory of $\der(A)$ containing $A$
and stable under forming direct factors. For each
$P$ in $\per(A)$ and each $M\in\der(A)$, we have canonical
isomorphisms
\[
D\RHom_A(P,M) \iso \RHom_A(M,D\RHom_A(P,A)) \mbox{ and }
P \lten_A DA \iso D\RHom_A(P,A).
\]
So we obtain a canonical isomorphism
\[
D\RHom_A(P,M) \iso \RHom_A(M, P\lten_A DA).
\]
Thus, if we are given a quasi-isomorphism of
dg $A$-$A$-bimodules
\[
\phi: A[n] \to DA \ko
\]
we obtain
\[
D\RHom_A(P,M) \iso \RHom_A(M, P[n])
\]
and in particular $\per(A)$ is Calabi-Yau of CY-dimension $n$.

{\bf (5)} To consider a natural application of the preceding example,
let $B$ be the symmetric algebra on a finite-dimensional vector space
$V$ of dimension $n$ and $\ct\subset \der(B)$ the localizing
subcategory generated by the trivial $B$-module $k$ (\ie the smallest
full triangulated subcategory stable under infinite sums and
containing the trivial module).  Let $\ct^c$ denote its subcategory of
compact objects.  This is exactly the triangulated subcategory of
$\der(B)$ generated by $k$, and also exactly the subcategory of the
complexes whose total homology is finite-dimensional and supported in
$0$.  Then $\ct^c$ is Calabi-Yau of CY-dimension $n$. Indeed, if
\[
A= \RHom_B(k,k)
\]
is the Koszul dual of $B$ (thus, $A$ is the exterior algebra on
the dual of $V$ concentrated in degree $1$; it is endowed with $d=0$),
then the functor
\[
\RHom_B(k,?): \der(B) \to \der(A)
\]
induces equivalences from $\ct$ to $\der(A)$ and $\ct^c$ to
$\per(A)$, \cf for example \cite{Keller94}. Now we have a canonical
isomorphism of $A$-$A$-bimodules $A[n] \iso DA$ so that
$\per(A)$ and $\ct$ are Calabi-Yau of CY-dimension $n$.
As pointed out by I.~Reiten, in this case, the Calabi-Yau
property even holds more generally:
Let $M\in\der(B)$ and denote by $M_\ct \to M$
the universal morphism from an object of $\ct$ to $M$.
Then, for $X\in\ct^c$, we have natural morphisms
\[
\Hom_{\der(B)}(M,X[n]) \to \Hom_{\ct}(M_\ct, X[n])
\iso D\Hom_\ct(X,M_\ct) \iso D\Hom_{\der(B)}(X,M).
\]
The composition
\[
(*) \quad\quad\quad \Hom_{\der(B)}(M,X[n]) \to D\Hom_{\der(B)}(X,M)
\]
is a morphism of (co-)homological functors in $X\in\ct^c$
(resp. $M\in\der(B)$). We claim that it is an isomorphism
for $M\in\per(B)$ and $X\in \ct^c$. It suffices to prove this for
$M=B$ and $X=k$. Then one checks it using the fact that
\[
\RHom_B(k,B) \iso k[-n].
\]
These arguments still work for certain non-commutative
algebras $B$: If $B$ is an Artin-Schelter regular algebra
\cite{ArtinSchelter87} \cite{ArtinTateVandenBergh90} of
global dimension $3$ and type $A$ and $\ct$ the localizing
subcategory of the derived category $\der(B)$ of non graded
$B$-modules generated by the trivial module, then $\ct^c$ is
Calabi-Yau and one even has the isomorphism $(*)$ for each perfect
complex of $B$-modules $M$ and each $X\in\ct^c$, \cf for example
section~12 of \cite{LuPalmieriWuZhang04}.

\subsection{Orbit categories with the Calabi-Yau property}
The main theorem yields the following

\begin{corollary} If $d\in \Z$ and $Q$ is a quiver whose underlying graph is
Dynkin of type $A$, $D$ or $E$, then
\[
\ct= \der^b(kQ)/\tau^{-1} S^{d-1}
\]
is Calabi-Yau of CY-dimension $d$.
In particular, the cluster category $\cc_{kQ}$ is Calabi-Yau
of dimension $2$ and the category of projective modules over
the preprojective algebra $\Lambda(Q)$ is Calabi-Yau of
CY-dimension $1$.
\end{corollary}

The category of projective modules over the preprojective algebra
$\Lambda(L_n)$ of example~\ref{ss:Projectives-over-Ln} does not fit
into this framework.  Nevertheless, it is also Calabi-Yau of
CY-dimension $1$, since we have $\tau=\id$ in this category and
therefore $\nu=S\tau=S$.

\subsection{Module categories over Calabi-Yau categories}
Calabi-Yau triangulated categories turn out to be
\lq self-reproducing\rq: Let $\ct$ be a triangulated category.
Then the category $\mod \ct$ of finitely generated functors
from $\ct\op$ to $\Mod k$ is abelian and Frobenius,
\cf \cite{Freyd66a}, \cite{Neeman99}. If we denote
by $\Sigma$ the exact functor $\mod\ct\to\mod\ct$
which takes $\Hom(?,X)$ to $\Hom(?,SX)$, then it is
not hard to show \cite{Freyd66a} \cite{Geiss04}
that we have
\[
\Sigma \iso S^3
\]
as triangle functors $\ul{\mod}\ct\to \ul{\mod}\ct$.
One deduces the following lemma, which is a variant
of a result which Auslander-Reiten \cite{AuslanderReiten96}
obtained using dualizing $R$-varieties \cite{AuslanderReiten74c}
and their functor categories \cite{AuslanderReiten73},
\cf also \cite{AuslanderReiten75a} \cite{Reiten02}.
A similar result is due to Geiss \cite{Geiss04}.

\begin{lemma} If $\ct$ is Calabi-Yau of CY-dimension $d$, then
the stable category $\ul{\mod}\,\ct$ is Calabi-Yau of
CY-dimension $3d-1$. Moreover, if the suspension of $\ct$
is of order $n$, the order of the suspension
functor of $\ul{\mod}\,\ct$ divides $3n$.
\end{lemma}

For example, if $A$ is the preprojective algebra of a
Dynkin quiver or equals $\Lambda(L_n)$, then we find that the stable category
$\ul{\mod}\,A$ is Calabi-Yau of CY-dimension $3\times 1 -1 = 2$.
This result, with essentially the same proof, is due
to Auslander-Reiten \cite{AuslanderReiten96}.
For the preprojective algebras of Dynkin quivers,
it also follows from a much finer result due to
Ringel and Schofield (unpublished). Indeed, they have proved
that there is an isomorphism
\[
\Omega_{A^e}^3(A) \iso DA
\]
in the stable category of bimodules, \cf Theorems 4.8 and
4.9 in \cite{BrennerButlerKing02}. This implies the
Calabi-Yau property since we also have an isomorphism
\[
DA\ten_A DA \iso A
\]
in the stable category of bimodules, by the remark following
definition 4.6 in \cite{BrennerButlerKing02}. For the
algebra $\Lambda(L_n)$, the analogous result follows
from Proposition 2.3 of \cite{BialkowskiErdmannSkowronski04}.

\section{Universal properties}
\label{s:Universal-properties}

\subsection{The homotopy category of small dg categories}
\label{ss:htpy-cat-dgcat}
Let $k$ be a field.
A {\em differential graded (=dg) $k$-module} is a
$\Z$-graded vector space
\[
V=\bigoplus_{p\in\Z} V^p
\]
endowed with a differential $d$ of degree $1$.
The {\em tensor product $V\ten W$}
of two dg $k$-modules is the graded space with components
\[
\bigoplus_{p+q=n} V^p\ten W^q \ko n\in\Z \ko
\]
and the differential $d\ten \id + \id\ten d$, where the tensor product
of maps is defined using the Koszul sign rule.  A {\em dg category}
\cite{Keller94} \cite{Drinfeld04} is a $k$-category $\ca$ whose
morphism spaces are dg $k$-modules and whose compositions
\[
\ca(Y,Z)\ten \ca(X,Y) \to \ca(X,Z)
\]
are morphisms of dg $k$-modules.  For a dg category $\ca$, the {\em
category $H^0(\ca)$} has the same objects as $\ca$ and has morphism
spaces $H^0\ca\,(X,Y)$, $X,Y\in\ca$.  A {\em dg functor $F:\ca\to\cb$}
between dg categories is a functor compatible with the grading and the
differential on the morphism spaces. It is a {\em quasi-equivalence}
if it induces quasi-isomorphisms in the morphism spaces and an
equivalence of categories from $H^0(\ca)$ to $H^0(\cb)$. We denote by
$\dgcat$ the category of small dg categories.  The {\em homotopy category}
of small dg categories is the localization $\Ho(\dgcat)$ of $\dgcat$
with respect to the class of quasi-equivalences. According to
\cite{Tabuada04}, the category $\dgcat$ admits a structure of Quillen
model category (\cf \cite{DwyerSpalinski95}, \cite{Hovey99}) whose
weak equivalences are the quasi-equivalences. This implies in
particular that for $\ca,\cb\in\dgcat$, the morphisms from $\ca$ to
$\cb$ in the localization $\Ho(\dgcat)$ form a set.

\subsection{The bimodule bicategory} For two dg categories
$\ca$, $\cb$, we denote by $\rep(\ca,\cb)$ the full subcategory
of the derived category $\der(\ca\op\ten\cb)$, \cf \cite{Keller94},
whose objects are the dg $\ca$-$\cb$-bimodules $X$ such that
$X(?,A)$ is isomorphic to a representable functor in $\der(\cb)$
for each object $A$ of $\ca$. We think of the objects of $\rep(\ca,\cb)$
as `representations up to homotopy' of $\ca$ in $\cb$.
The {\em bimodule bicategory $\bim$},
\cf \cite{Keller94} \cite{Drinfeld04},
has as objects all small dg categories; the morphism category between two
objects $\ca$, $\cb$ is $\rep(\ca,\cb)$; the composition
bifunctor
\[
\rep(\cb,\cc)\times \rep(\ca,\cb) \to \rep(\ca,\cb)
\]
is given by the derived tensor product $(X,Y) \mapsto X\lten_\cb Y$.
For each dg functor $F:\ca\to\cb$,
we have the dg bimodule
\[
X_F : A \mapsto \cb\,(?,FA) \ko
\]
which clearly belongs to $\rep(\ca,\cb)$. One can show that
the map $F \mapsto X_F$ induces a bijection, compatible
with compositions, from the
set of morphisms from $\ca$ to $\cb$ in $\Ho(\dgcat)$ to the
set of isomorphism classes of bimodules $X$ in $\rep(\ca,\cb)$.
In fact, a much stronger result by B.~To\"en \cite{Toen04} relates
$\bim$ to the Dwyer-Kan localization of $\dgcat$.

\subsection{Dg orbit categories} \label{ss:dg-orbit-categories}
Let $\ca$ be a small dg category and $F\in\rep(\ca,\ca)$.
We assume, as we may, that $F$ is given by a cofibrant bimodule.
For a dg category $\cb$, define $\tilde{\eff}_0(\ca, F,\cb)$
to be the category whose objects are the pairs formed by
an $\ca$-$\cb$-bimodule $P$ in $\rep(\ca,\cb)$ and
a morphism of dg bimodules
\[
\phi: P \to PF.
\]
Morphisms are the morphisms of dg bimodules $f: P \to P'$
such that we have $\phi'\circ f = (fF)\circ \phi$ in the
category of dg bimodules. Define $\eff_0(\ca, F,\cb)$
to be the localization of $\tilde{\eff}_0(\ca,F,\cb)$ with
respect to the morphisms $f$ which are quasi-isomorphisms of
dg bimodules. Denote by $\eff(\ca,F,\cb)$
the full subcategory of $\eff_0(\ca,F,\cb)$ whose objects are the
$(P,\phi)$ where $\phi$ is a quasi-isomorphism.
It is not hard to see that the assignments
\[
\cb\mapsto \eff_0(\ca, F,\cb) \quad \mbox{and} \quad
\cb\mapsto \eff(\ca, F,\cb)
\]
are $2$-functors from $\rep$ to the category of small
categories.

\begin{theorem}
\begin{itemize}
\item[a)] The $2$-functor $\eff_0(\ca,F,?)$ is $2$-representable,
\ie there is small dg category $\cb_0$ and a pair $(P_0,\phi_0)$
in $\eff(\ca,F,\cb_0)$ such that for each small dg category $\cb$,
the functor
\[
\rep(\cb_0,\cb) \to \eff_0(\ca,F,\cb_0) \ko G \mapsto G\circ P_0
\]
is an equivalence.
\item[b)] The $2$-functor $\eff(\ca,F,?)$ is $2$-representable.
\item[c)] For a dg category $\cb$, a pair $(P,\phi)$ is a
$2$-representative for $\eff_0(\ca, F, ?)$ iff $H^0(P): H^0(\ca) \to H^0(\cb)$
is essentially surjective and, for all objects $A,B$ of $\ca$,
the canonical morphism
\[
\bigoplus_{n\in\N} \ca\,(F^n A, B) \to \cb\,(PA, PB)
\]
is invertible in $\der(k)$.
\item[d)] For a dg category $\cb$, a pair $(P,\phi)$ is a
$2$-representative for $\eff(\ca,F,?)$ iff $H^0(P): H^0(\ca) \to H^0(\cb)$
is essentially surjective and, for all objects $A,B$ of $\ca$, 
the canonical morphism
\[
\bigoplus_{c\in\Z}\colim_{r\gg 0}
\ca\,(F^{p+r} A,F^{p+c+r} B) \to \cb\,(PA, PB)
\]
is invertible in $\der(k)$.
\end{itemize}
\end{theorem}

We {\em define $\ca/F$} to be the $2$-representative of
$\eff(\ca,F,?)$. For example, in the notations of~\ref{ss:dgorbitcategory},
$\ca/F$ is the dg orbit category $\cb$. It follows from part d) of
the theorem that we have an equivalence
\[
H^0 (\ca)/ H^0(F) \to H^0 (\ca/F).
\]

\begin{proof} We only sketch a proof and refer to \cite{Tabuada06}
for a detailed treatment. Define $\cb_0$ to be the dg category with
the same objects as $\ca$ and with the morphism spaces
\[
\cb_0\,(A,B) = \bigoplus_{n\in\N} \ca\,(F^n A, B).
\]
We have an obvious dg functor $P_0:\ca\to\cb_0$ and an obvious
morphism $\phi : P_0  \to P_0 F$. The pair $(P_0, \phi_0)$ is
then $2$-universal in $\rep$. This yields a) and c). For b), one
adjoins a formal homotopy inverse of $\phi$ to $\cb_0$. One obtains
d) by computing the homology of the morphism spaces in the resulting
dg category.
\end{proof}

\subsection{Functoriality in $(\ca,F)$.} Let a square of
$\rep$
\[
\xymatrix{\ca \ar[r]^G \ar[d]_F & \ca' \ar[d]^{F'} \\
\ca \ar[r]^G & \ca'}
\]
be given and an isomorphism
\[
\gamma : F'G \to GF
\]
of $\rep(\ca,\ca')$. We assume, as we may, that $\ca$ and $\ca'$ are
cofibrant in $\dgcat$ and that $F$, $F'$ and $G$ are given
by cofibrant bimodules. Then $F'G$ is a cofibrant bimodule
and so $\gamma: F'G \to GF$ lifts to a morphism of
bimodules
\[
\tilde{\gamma}: F'G \to GF.
\]
If $\cb$ is another dg category and $(P,\phi)$ an
object of $\eff_0(\ca', F', \cb)$, then the composition
\[
\xymatrix{PG \ar[r]^{\phi G} & PF'G \ar[r]^{P\tilde{\gamma}} & PGF}
\]
yields an object $PG \to PGF$ of $\eff_0(\ca,F,\cb)$. Clearly,
this assignment extends to a functor, which induces a functor
\[
\eff(\ca',F',\cb) \to \eff(\ca,F,\cb).
\]
By the 2-universal property of section~\ref{ss:dg-orbit-categories},
we obtain an induced morphism
\[
\overline{G}: \ca/F \to \ca'/F'.
\]
One checks that the composition of two pairs $(G,\gamma)$ and
$(G',\gamma')$ induces a functor isomorphic to the composition
of $\overline{G}$ with $\overline{G'}$.

\subsection{The bicategory of enhanced triangulated categories}
We refer to \cite[2.1]{Keller99} for the notion of an
{\em exact dg category}. We also call these categories {\em
pretriangulated} since if $\ca$ is an exact dg category, then
$H^0(\ca)$ is triangulated. More precisely, $\ce=Z^0(\ca)$ is
a Frobenius category and $H^0(\ca)$ is its associated stable
category $\ul{\ce}$ (\cf example (3) of section~\ref{ss:CY-examples}
for these notions).

The inclusion of the full subcategory of (small) exact dg categories
into $\Ho(\dgcat)$ admits a left adjoint, namely the functor
$\ca \mapsto \pretr(\ca)$ which maps a dg category to its
`pretriangulated hull' defined in \cite{BondalKapranov91}, \cf
also \cite[2.2]{Keller99}. More precisely,
the adjunction morphism $\ca\to\pretr(\ca)$ induces
an equivalence of categories
\[
\rep\,(\pretr(\ca),\cb) \to \rep(\ca,\cb)
\]
for each exact dg category $\cb$, \cf \cite{Tabuada05}.

The {\em bicategory $\enh$ of enhanced \cite{BondalKapranov91}
triangulated categories}, \cf \cite{Keller94} \cite{Drinfeld04},
has as objects all small exact dg categories; the morphism category
between two objects $\ca$, $\cb$ is $\rep(\ca,\cb)$; the composition
bifunctor
\[
\rep(\cb,\cc)\times \rep(\ca,\cb) \to \rep(\ca,\cc)
\]
is given by the derived tensor product $(X,Y) \mapsto X\lten_\cb Y$.

\subsection{Exact dg orbit categories}
\label{ss:exact-dg-orbit-categories}
Now let $\ca$ be an exact dg category and $F\in\rep(\ca,\ca)$.
Then $\ca/F$ is the dg orbit category of subsection~\ref{ss:dgorbitcategory}
and $\pretr(\ca/F)$ is an exact dg category such that
$H^0\pretr(\ca/F)$ is the triangulated hull of
section~\ref{s:triangulated-hull}. In particular, we obtain that
the triangulated hull is the stable category of a Frobenius
category. From the construction, we
obtain the universal property:

\begin{theorem}
For each exact dg category $\cb$, we have an equivalence
of categories
\[
\rep(\pretr(\ca/F), \cb) \to \eff(F,\cb).
\]
\end{theorem}

\subsection{An example} Let $A$ be a finite-dimensional
algebra of finite global dimension and $TA$ the trivial extension algebra,
\ie the vector space $A\oplus DA$ endowed with the multiplication
defined by
\[
(a,f)(b,g)=(ab,ag+fb)\ko (a,f),(b,g)\in TA\ko
\]
and the grading such that $A$ is in degree $0$ and $DA$
in degree $1$. Let $F: \der^b(A) \to \der^b(A)$ equal $\tau S^2=\nu S$ and
let $\tilde{F}$ be the dg lift of $F$ given by $?\ten_A R[1]$, where
$R$ is a projective bimodule resolution of $DA$. Let $\der^b(A)_{dg}$
denote a dg category quasi-equivalent to the dg category of
bounded complexes of finitely generated projective $A$-modules.

\begin{theorem}
The following are equivalent
\begin{itemize}
\item[(i)] The $k$-category $\der^b(A)/F$ is naturally equivalent
to its \lq triangulated hull\rq
\[
H^0(\pretr (\der^b(A)_{dg}/\tilde{F})) .
\]
\item[(ii)] Each finite-dimensional $TA$-module admits a grading.
\end{itemize}
\end{theorem}

\begin{proof}
We have a natural functor
\[
\mod A \to \grmod TA
\]
given by viewing an $A$-module as a graded $TA$-module
concentrated in degree~$0$. As shown by D.~Happel \cite{Happel87},
\cf also \cite{KellerVossieck87},
this functor extends to a triangle equivalence $\Phi$
from $\der^b(A)$ to the stable category $\ul{\grmod}\, TA$,
obtained from $\grmod TA$ by killing all morphisms factoring
through projective-injectives. We would like to show that
we have an isomorphism of triangle functors
\[
\Phi \circ \tau S^2 \iso \Sigma \circ \Phi
\]
where $\Sigma$ is the grading shift functor for graded
$TA$-modules: $(\Sigma M)^p=M^{p+1}$ for all $p\in\Z$.
From \cite{Happel87}, we know that $\tau S \iso \nu$, where
$\nu = ?\lten_A DA$. Thus it remains to show that
\[
\Phi \circ \nu S \iso \Sigma\circ \Phi.
\]
As shown in \cite{KellerVossieck87}, the equivalence $\Phi$
is given as the composition
\begin{align*}
\der^b(\mod A) \to \der^b & (\grmod TA)  \to \\
 & \der^b(\grmod TA)/\per(\grmod TA) \to \ul{\grmod}\,TA \ko
\end{align*}
where the first functor is induced by the above inclusion,
the notation $\per(\grmod TA)$ denotes the triangulated subcategory
generated by the projective-injective $TA$-modules and
the last functor is the \lq stabilization functor\rq\, \cf
\cite{KellerVossieck87}. We have a short exact sequence of
graded $TA$-modules
\[
0 \to \Sigma^{-1}(DA) \to TA \to A \to 0.
\]
We can also view it as a sequence of left $A$ and right
graded $TA$-modules. Let $P$ be a bounded complex of projective
$A$-modules. Then we obtain a short exact sequence
of complexes of graded $TA$-modules
\[
0 \to \Sigma^{-1}(P\ten_A DA) \to P\ten_A TA \to P \to 0
\]
functorial in $P$. It yields a functorial triangle in
$\der^b(\grmod A)$. The second term belongs to $\per(\grmod TA)$.
Thus in the quotient category
\[
\der^b(\grmod TA)/\per(\grmod TA) \ko
\]
the triangle reduces to a functorial isomorphism
\[
P \iso S \Sigma^{-1} \nu P.
\]
Thus we have a functorial isomorphism
\[
\Phi(P) \iso S \Sigma^{-1} \Phi(\nu P).
\]
Since $A$ is of finite global dimension, $\der^b(\mod A)$ is
equivalent to the homotopy category of bounded complexes of
finitely generated projective $A$-modules. Thus we get the required
isomorphism
\[
\Sigma \Phi \iso \Phi S \nu.
\]
More precisely, one can show that $\underline{\grmod}\, TA$ has a canonical
dg structure and that there is an isomorphism
\[
(\der^b(A))_{dg} \iso (\ul{\grmod}\, TA)_{dg}
\]
in the homotopy category of small dg categories which
induces Happel's equivalence and under which $\Sigma$
corresponds to the lift $\tilde{F}$ of $F=S\nu$.
Hence the orbit categories $\der^b(\mod A)/\tau S^2$ and
$\ul{\grmod}\, TA/\Sigma$ are equivalent and we
are reduced to determining when $\ul{\grmod}\, TA/\Sigma$
is naturally equivalent to its triangulated hull.
Clearly, we have a full embedding
\[
\ul{\grmod}\, TA/\Sigma \to \ul{\mod}\, TA
\]
and its image is formed by the $TA$-modules which admit a
grading. Now $\ul{\mod}\,TA$ is naturally equivalent to
the triangulated hull. Therefore, condition (i) holds
iff the embedding is an equivalence iff each finite-dimensional
$TA$-module admits a grading.
\end{proof}

In \cite{Skowronski87}, A.~Skowro\'nski has produced a class of examples where
condition (ii) does not hold. The simplest of these is the algebra
$A$ given by the quiver
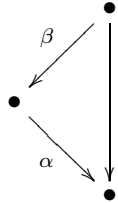
\[
\xymatrix{                           & \bullet \ar[dl]_{\beta} \ar[dd] \\
            \bullet \ar[dr]_{\alpha} &  \\
                                     & \bullet
}
\]
with the relation $\alpha\beta=0$. Note that this algebra is
of global dimension~$2$.

\subsection{Exact categories and standard functors} Let $\ce$ be a small exact
$k$-category. Denote by $\cc^b(\ce)$ the category of bounded complexes
over $\ce$ and by $\ac^b(\ce)$ its full subcategory formed by
the acyclic bounded complexes. The categories with the
same objects but whose morphisms are given by the
morphism {\em complexes} are denoted respectively by
$\cc^b(\ce)_{dg}$ and $\ac^b(\ce)_{dg}$. They are
exact dg categories and so is the dg quotient
\cite{Keller99} \cite{Drinfeld04}
\[
\der^b(\ce)_{dg} = \cc^b(\ce)_{dg}/\ac^b(\ce)_{dg}.
\]
Let $\ce'$ be another small exact $k$-category.
We call a triangle functor $F:\der^b(\ce) \to \der^b(\ce')$ a
{\em standard functor} if it is isomorphic to the triangle
functor induced by a morphism
\[
\tilde{F} : \der^b(\ce)_{dg} \to \der^b(\ce')_{dg}
\]
of $\Ho(\dgcat)$. Slightly abusively,
we then call $\tilde{F}$ a {\em dg lift} of $F$.
Each exact functor $\ce \to \ce'$ yields a standard functor;
a triangle functor is standard iff it admits a lift to an object of
$\rep(\der^b(\ce)_{dg}, \der^b(\ce')_{dg})$;
compositions of standard functors are standard;
an adjoint (and in particular, the inverse) of a standard functor is
standard.

If $F: \der^b(\ce) \to \der^b(\ce)$ is
a standard functor with dg lift $\tilde{F}$, we have the dg orbit
category $\der^b(\ce)_{dg}/\tilde{F}$ and its pretriangulated hull
\[
\der^b(\ce)_{dg}/\tilde{F} \to \pretr(\der^b(\ce)_{dg}/\tilde{F}).
\]
The examples in section~\ref{s:examples} show that
this functor is not an equivalence in general.

\subsection{Hereditary categories}
\label{ss:hereditary-categories} Now suppose that $\ch$ is
a small hereditary abelian $k$-category with the Krull-Schmidt
property (indecomposables have local endomorphism rings and
each object is a finite direct sum of indecomposables)
where all morphism and extension spaces are
finite-dimensional. Let
\[
F: \der^b(\ch) \to \der^b(\ch)
\]
be a standard functor with dg lift $\tilde{F}$.

\begin{theorem} Suppose that $F$ satisfies assumptions
2) and 3) of the main theorem in section~\ref{s:main-theorem}.
Then the canonical functor
\[
\der^b(\ch)/F
\to
H^0(\pretr(\der^b(\ch)_{dg}/\tilde{F}))
\]
is an equivalence of $k$-categories. In particular, the orbit
category $\der^b(\ch)/F$ admits a triangulated structure
such that the projection functor becomes a triangle functor.
\end{theorem}

The proof is an adaptation, left to the reader, of the proof of
the main theorem.

Suppose for example that $\der^b(\ch)$ has a Serre functor $\nu$. Then
$\nu$ is a standard functor since it is induced by the tensor
product with bimodule
\[
(A,B) \mapsto D\Hom_{\der^b(\ch)_{dg}}(B,A) \ko
\]
where $D=\Hom_k(?,k)$. The functor $\tau^{-1} = S\nu^{-1}$ induces
equivalences $\ci \iso S\cp$ and $\ch_{ni} \to \ch_{np}$, where
$\cp$ is the subcategory of projectives, $\ci$ the subcategory
of injectives, $\ch_{np}$ the subcategory of objects without
a projective direct summand and $\ch_{ni}$ the subcategory of
objects without an injective direct summand. Now let $n\geq 2$
and consider the autoequivalence
$F=S^{n} \nu^{-1}=S^{n-1} \tau^{-1}$ of $\der^b(\ch)$. Clearly
$F$ is standard. It is not hard to see that $F$ satisfies the
hypotheses 2) and 3) of the main theorem in section~\ref{s:main-theorem}.
Thus the orbit category
\[
\der^b(\ch)/F = \der^b(\ch)/S^{n} \nu^{-1}
\]
is triangulated. Note that we have excluded the case $n=1$ since
the hypotheses 2) and 3) are not satisfied in this case, in general.


\def\cprime{$'$}
\providecommand{\bysame}{\leavevmode\hbox to3em{\hrulefill}\thinspace}
\providecommand{\MR}{\relax\ifhmode\unskip\space\fi MR }
\providecommand{\MRhref}[2]{%
  \href{http://www.ams.org/mathscinet-getitem?mr=#1}{#2}
}
\providecommand{\href}[2]{#2}

\end{document}